\documentclass[12pt]{article}

\usepackage{amssymb,amsmath,theorem}
\usepackage{epsfig}
\usepackage{calc,ifthen}
\usepackage{enumerate}
\usepackage{url}

\reversemarginpar 

\newcommand\al{\alpha}
\renewcommand\th{\theta} 
\newcommand\fld{{\mathbb F}}
\newcommand\cA{{\mathcal A}}

\newcommand\sbs{\subseteq}
\newcommand\seq[3]{#1_{#2},\ldots,#1_{#3}}
\newcommand\comp[1]{{\mkern2mu\overline{\mkern-2mu#1}}}

\newcommand\tr{\operatorname{tr}}
\newcommand\summ{\operatorname{sum}}

\theoremstyle{change}
{\theorembodyfont{\slshape}
\newtheorem{theorem}{Theorem.}[section]
\newtheorem{lemma}[theorem]{Lemma.}
\newtheorem{corollary}[theorem]{Corollary.}}
\theorembodyfont{\rmfamily}

\def\sqr#1#2{{\vbox{\hrule height.#2pt
     \hbox{\vrule width.#2pt height#1pt \kern#1pt
         \vrule width.#2pt}\hrule height.#2pt}}}
\def\eqed{\sqr53}
\def\qed{%
     \ifmmode\eqno\eqed
     \else\nobreak\ \hfill\eqed\medbreak\fi}

\newcommand\Cref[1]{Chapter~\ref{cha:#1}}

\newcommand\lref[1]{Lemma~\ref{lem:#1}}
\newcommand\tref[1]{Theorem~\ref{thm:#1}}
\newcommand\cref[1]{Corollary~\ref{cor:#1}}
\newcommand\sref[1]{Section~\ref{sec:#1}}

\newcommand{\diag}{\operatorname{diag}}
\newcommand{\psd}{\succeq0}
\newcommand{\norm}[1]{\left\lVert#1\right\rVert}

\newcommand{\pp}[1]{\left(#1\right)}

\newcommand{\cc}[1]{\left\{#1\right\}}

\newcommand{\mat}[1]{\left(\begin{matrix} #1 \end{matrix}\right)}

\newcommand{\1}{\mathbf{1}}

\title{Eigenvalue bounds for independent sets}

\author{
  C. D. Godsil$\mbox{}^1$\footnotemark[1]~ and M. W. Newman$\mbox{}^2$\footnotemark[1]\footnote{Research supported by NSERC.}\\
  \normalsize $\mbox{}^1$Department of Combinatorics and Optimization\\
  \normalsize University of Waterloo, CANADA\\
  \normalsize $\mbox{}^2$School of Mathematical Sciences\\
  \normalsize Queen Mary, University of London, UK
}

\begin{document}

\maketitle

\abstract{We derive bounds on the size of an independent set based on eigenvalues. This generalizes a result
due to Delsarte and Hoffman. We use this to obtain new bounds on the independence number of the
Erd\H{o}s-R\'{e}nyi graphs. We investigate further properties of our bounds, and show how our results on the
Erd\H{o}s-R\'{e}nyi graphs can be extended to other polarity graphs.}

\section{Introduction}

Let $\fld$ be a finite field of order $q$ and let $V$ be a 3-dimensional vector space over $\fld$.  The
1-dimensional subspaces of $V$ are the points of the projective plane $PG(2,q)$, and the 2-dimensional
subspaces are the lines.  It follows that each point can be represented by a non-zero vector, namely any
vector that spans the corresponding 1-dimensional subspace.  Two points $a$ and $b$, represented by vectors
$x$ and $y$ respectively, are \textsl{orthogonal} if $x^Ty=0$.  The Erd\H{o}s-R\'enyi graph $ER(q)$ is the
graph with the points of $PG(2,q)$ as its vertices, where two vertices are adjacent if and only if they are
orthogonal.

The graph $ER(q)$ has $q^2+q+1$ vertices and each vertex has exactly $q+1$ neighbours.  There is one problem
though: by standard results in finite geometry, there are exactly $q+1$ vertices that are adjacent to
themselves.  Thus our graph has $q+1$ loops.  The Erd\H{o}s-R\'enyi graphs are of interest because they do not
contain any 4-cycles, but nonetheless they have a large number of edges; this is the motivation for
\cite{ERS66}. For further work on these graphs, see \cite{Par76,Williford}.

Our object in this paper is to derive good bounds on the size of an independent set.  (Where an independent
set is a subset of the vertices such that no two distinct vertices are adjacent.  Thus an independent set may
contains vertices with loops.)  Let $\al(X)$ denote the maximum number of vertices in an independent set in
$X$.  There is a standard bound for $\al(X)$ in terms of the eigenvalues of $X$, due to Delsarte and Hoffman
(see~\cite[Section~3.3]{Delsarte} or~\cite[Page 115]{CDS}; alternatively~\cite{mike} for more recent
work). However this bound only applies to regular graphs with no loops, and consequently our first task in
this paper is to derive an extension of it.  With this in hand we are able to derive new (and better) bounds
on the size of independent sets in $ER(q)$.

We conclude the paper by describing some more general classes of graphs to which our new bound can be applied.
These are obtained as follows.  Suppose $Y$ is a connected $k$-regular bipartite graph on $2v$ vertices, and
suppose there is an automorphism $\th$ of $Y$ with order two that swaps the two colour classes of $Y$.  The
\textsl{quotient graph} $Y/\th$ is a graph with the $v$ orbits of $\th$ as its vertices, and with $b_{i,j}$
arcs from orbit $i$ to orbit $j$, where $b_{i,j}$ is the number of edges in $Y$ from a vertex in orbit $i$ to
the vertices of orbit $j$.  Any orbit of $\th$ that contains two adjacent vertices gives rise to a vertex in
$Y/\th$ with a loop.  It is not hard to show that $b_{j,i}=b_{i,j}$ and that, if $Y$ has no 4-cycles, then
$Y/\th$ does not have any multiple edges.  The graph $ER(q)$ can be constructed in this way from the incidence
graph of the projective plane $PG(2,q)$; the vertices of the incidence graph are the points and lines of
$PG(2,q)$ and a point is adjacent to a line in the graph if it is incident with the line in the geometry.  The
map that sends the point represented by the non-zero vector $x$ to the line consisting of the points
represented by the non-zero vectors $y$ such that $x^Ty=0$ gives rise to an automorphism $\th$ of order two
that swaps points and lines.

\section{General Framework}\label{sec:framework}

We concern ourselves with bounding the size of an independent set in a graph. We will permit loops on
vertices, but we will allow these vertices to be included in an independent set; in other words, we define an
independent set to be a set of vertices of which no two distinct members are adjacent. Allowing loops is the
more general option. If we later wish to exclude them we may delete the looped vertices.

We will need a little linear algebra. Recall that a symmetric matrix $B$ is \emph{positive semidefinite} if
all of its eigenvalues are non-negative; equivalently, if $x^TBx \geq 0$ for all vectors $x$. We write $B
\psd$. If $B$ is positive semi-definite, then $x^TBx=0$ if and only if $Bx=0$.

Let $X$ be a graph with vertex set $V$, $|V|=n$, possibly containing loops, and let $A$ be its adjacency
matrix. Let $T=\diag\pp{\seq{t}{1}{n}}$ be a diagonal matrix such that $T+A\psd$. Also, let $d_i$ be the
degree of vertex $i$, with loops counted once each. Consider an independent set $S$ of size $s$, and let $s_1$
be the number of loops on vertices in $S$. Let $z$ be the characteristic vector of $S$. Then we have
\[	\pp{z-\frac{s}{n}\1}^T \pp{T+A} \pp{z-\frac{s}{n}\1} \geq 0 \]
Expanding this we obtain the following result.
\begin{lemma}\label{lem:genbound}
  Let $X$ be a graph with vertex set $V$, $|V|=n$ and vertex degrees $\seq{d}{1}{n}$. Let $A$ be its adjacency
  matrix, and let $T=\diag\pp{\seq{t}{1}{n}}$ be such that $T+A\psd$. If $S$ is an independent set of size $s$
  containing $s_1$ loops, then:
  \[
  \frac{s^2}{n^2} \sum_{i \in V} \pp{t_i+d_i}
  - 2\frac{s}{n} \sum_{i \in S} \pp{t_i+d_i} + \sum_{i \in S} t_i \geq -s_1.\qed
  \]
\end{lemma}
This gives a bound on $s$. However, it is difficult to apply in general, partly because the sums depend not
only on $s$ but on $S$. Furthermore, the bound obtained will depend on the choice of $T$. We do not know how
to choose $T$ optimally (or even if there is a single optimal choice for all graphs). In the present paper, we
consider specific choices for $T$.

It will be useful to define the following parameters of a set $S$:
\begin{align}
  \comp{d}_S &= \frac{1}{s}\sum_{i \in S} d_i, \notag\\
  k_S &= 2 \comp{d}_S - \frac{1}{n}\sum_{i \in V} d_i. \notag
\end{align}
Note that for $k$-regular graphs, $\comp{d}_S = k_S = k$. It will be seen that these two parameters behave, in
some circumstances, as analogues to the degree of a regular graph.

\section{Bounds}\label{sec:bounds}

We consider two particular choices for $T$, producing bounds which we can regard as coming from, respectively,
the adjacency matrix and the Laplacian matrix of the graph.

\subsection{Adjacency Matrix}

If we let $\tau$ be the least eigenvalue of $A$, then we may set $T=-\tau I$ giving $T+A = A-\tau I
\psd$. If $X$ is regular and loopless, then simplification of \lref{genbound} gives the Delsarte-Hoffman bound.
\begin{corollary}\label{cor:rbound}
  Let $X$ be a $k$-regular graph with no loops, and $\tau$ the least eigenvalue of its adjacency matrix. For
  any independent set $S$ of size $s$, we have:
  \[ s \leq n \frac{-\tau}{k-\tau}.\qed \]
\end{corollary}

It turns out that $X$ need be neither regular nor loopless.

\begin{corollary}\label{cor:Abound}
  Let $X$ be a graph with no loops, and $\tau$ the least eigenvalue of its adjacency matrix. For any
  independent set $S$ of size $s$, we have:
  \[ s \leq n \frac{-\tau}{k_S-\tau}.\qed \]
\end{corollary}
Note that $k_S$ plays an analogous role to that of the degree in \cref{rbound}. However, $k_S$ can be zero or
even negative: the bound is then useless. To be precise, one should say that \cref{Abound} does not bound the
size of an independent set, but provides a family of bounds, one for each value of $k_S$ (or equivalently, one
bound for each value of $\comp{d}_S$). Any lower bound on $\comp{d}_S$ (such as the minimum degree) can be
used to make \cref{Abound} into a true bound on $s$.

If $X$ has loops, then \lref{genbound} is a nontrivial quadratic. The bounds we get are slightly messier, but
again, in the non-regular case, $k_S$ plays a role analogous to that of the degree in the regular case.

\begin{corollary}\label{cor:Abound1}
  Let $X$ be a $k$-regular graph with loops, and $\tau$ the least eigenvalue of its adjacency matrix. For any
  independent set $S$ of size $s$ containing $s_1$ loops, we have:
  \[ s \leq n \frac{ -\tau +\sqrt{\tau^2 + 4s_1\frac{k-\tau}{n}} }{ 2\pp{k-\tau} }.\qed \]
\end{corollary}

\begin{corollary}\label{cor:Abound2}
  Let $X$ be a graph with loops, and $\tau$ the least eigenvalue of its adjacency matrix. For any independent
  set $S$ of size $s$ containing $s_1$ loops , we have:
  \[ s \leq n \frac{ -\tau +\sqrt{\tau^2 + 4s_1\frac{k_S-\tau}{n}} }{ 2\pp{k_S-\tau} }.\qed \]
\end{corollary}

\subsection{Laplacian Matrix}\label{sec:Lapbounds}

For a graph $X$ with adjacency matrix $A$ and diagonal matrix of degrees $D$, recall that $L=D-A$ is the
Laplacian matrix of $X$. We always have $L \psd$, and in fact $0$ is an eigenvalue of multiplicity equal to
the number of components of $X$. The greatest eigenvalue of $L$ is at most twice the maximum degree; it is
also bounded by the number of vertices (see for instance~\cite{AM85}). If $X$ is regular then $L=kI-A$ and the
eigenvalues of $L$ and $A$ contain the same information. Accordingly we expect to recover previous bounds for
regular graphs and hope to obtain new ones in the non-regular case. Note that graphs that differ only by the
presence or absence of loops have the same Laplacian matrix. Thus, without loss of generality, we can assume
that the graph has no loops and set $s_1=0$.

If we let $\mu$ be the greatest eigenvalue of $L$, then we may set $T=\mu I -D$ giving $T+A = \mu I - L
\psd$. If the graph is regular, then we recover \cref{rbound}, as expected. If it is not regular, then we
obtain the following bound.
\begin{corollary}\label{cor:Lbound}
  Let $X$ be a loopless graph, and $\mu$ the greatest eigenvalue of its Laplacian matrix. For any independent set
  $S$ of size $s$, we have:
  \[ s \leq n \frac{\mu-\comp{d}_S}{\mu}.\qed \]
\end{corollary}
Note that $\comp{d}_S$ plays an analogous role to that of the degree in \cref{rbound}. We could have stated
\cref{Lbound} for graphs with loops, but it is more convenient to leave them out.

This result generalizes \cref{rbound}, but it also generalizes \cref{Abound1}. Let $X$ be a $k$-regular graph
with loops and let $Y$ be $X$ with the loop-edges deleted. The Laplacian eigenvalues of $X$ and $Y$ are
identical, so
\[ \mu(X) = \mu(Y) = k-\tau(Y). \]
Also,
\[ \comp{d}_S(X) = k - \frac{s_1}{s}. \]
Substituting into \cref{Lbound} yields \cref{Abound1}.

We can weaken \cref{Lbound} slightly to a more usable form, by noting that $\comp{d}_S \geq \delta$.
\begin{corollary}\label{cor:Lbound2}
  Let $X$ be any graph with minimum degree $\delta$, and $\mu$ the greatest eigenvalue of its Laplacian
  matrix. For any independent set $S$ of size $s$, we have:
  \[ s \leq n \frac{\mu-\delta}{\mu}.\qed \]
\end{corollary}

\section{Equality}\label{sec:equality}

If \lref{genbound} holds with equality then it follows that
\[ (T+A)(z-\frac{s}{n}\1) = 0, \]
and we have an eigenvector for $T+A$. Unpacking this equality gives a proof of the following.
\begin{lemma}\label{lem:gentight}
  Let $X$ be a graph with vertex degrees $\seq{d}{1}{n}$ and adjacency matrix $A$. Let
  $T=\diag(\seq{t}{1}{n})$ be such that $T+A\psd$. If $S$ is a set of $s$ vertices with no two distinct
  vertices adjacent such that \lref{genbound} holds with equality, then:
  \begin{enumerate}[(a)]
    \item Each vertex $i$ in $S$ has degree $d_i = t_i(\frac{n}{s}-1)$.
      \item Each vertex $i$ not in $S$ has $(d_i + t_i)\tfrac{s}{n}$ neighbours in $S$.\qed
  \end{enumerate}
\end{lemma}

\subsection{Adjacency matrix}

If $S$ is an independent set in a $k$-regular graph such that \cref{rbound} holds with equality, then we have
an eigenvector for $A$.
\[ A \pp{z-\frac{s}{n}\1} = \tau \pp{z-\frac{s}{n}\1}. \]
Since the vector $\1$ spans the $k$-eigenspace, $z$ lies in the sum of the greatest eigenspace
and the least eigenspace.

Using \lref{gentight}, we see that the bipartite subgraph induced by the partition $\{S,V(X)\setminus S\}$ is
semi-regular: vertices in $S$ have $-\tau(\frac{n}{s}-1)$ neighbours not in $S$, and vertices not in $S$ have
$(k-\tau)\frac{s}{n}$ neighbours in $S$. But vertices in $S$ have $k$ neighbours not in $S$, since the graph
is $k$-regular. It follows that vertices not in $S$ have $-\tau$ neighbours in $S$. The partition
$\{S,V(X)\setminus S\}$ is in fact equitable.

We summarize our findings in the following, due to Delsarte and Hoffman (unpublished).
\begin{theorem}\label{thm:Aboundeq}
  Let $X$ be a $k$-regular graph with no loops, and $\tau$ the least eigenvalue of its adjacency matrix. For
  any independent set $S$ of size $s$ and characteristic vector $z$, we have:
  \[ s \leq n \frac{-\tau}{k-\tau}.\]
  Furthermore, the following are equivalent:
  \begin{enumerate}[(a)]
    \setlength{\itemsep}{0pt}
    \setlength{\parskip}{0pt}
    \item Equality holds.
    \item $z$ is a linear combination of a $k$-eigenvector and a $\tau$-eigenvector.
    \item The bipartite subgraph induced by the partition $\{S,V(X)\setminus S\}$ is semi-regular.
    \item The partition $\{S, V(X) \setminus S\}$ is equitable.\qed
  \end{enumerate}
\end{theorem}

The $k$-eigenvectors of a regular graph are exactly the constant vectors. For non-regular graphs, the bounds
based on the adjacency matrix seem less useful. Specifically, if equality holds in \cref{Abound},
\cref{Abound1}, or \cref{Abound2}, then we still have an eigenvector for $A$, but as $\1$ is no longer an
eigenvector, the conditions of \lref{gentight} are not as useful.

\subsection{Laplacian matrix}

If $S$ is an independent set in a graph such that \cref{Lbound} holds with equality, then we have an
eigenvector for $L$:
\[ L \pp{z-\frac{s}{n}\1} = \mu \pp{z-\frac{s}{n}\1}. \]
Since the vector $\1$ spans the $0$-eigenspace, we have that $z$ lies in the sum of the least eigenspace and
the greatest eigenspace. Again, from \lref{gentight}, we see that the bipartite subgraph induced by
the partition $\{S,V(X)\setminus S\}$ is semi-regular: vertices in $S$ have $\mu(1-\frac{s}{n})$ neighbours
outside of $S$ and vertices outside of $S$ have $\mu\frac{s}{n}$ neighbours in $S$.

This does not quite say that $\{S,V(X) \setminus S\}$ is an equitable partition: the missing condition needed
is that every vertex in $V(X) \setminus S$ would have a constant number of neighbours in $V(X) \setminus
S$. But if this condition were to hold then the graph would be regular.

We summarize our findings as follows.
\begin{theorem}\label{thm:Lboundeq}
  Let $X$ be a graph with no loops, and $\mu$ the greatest eigenvalue of its Laplacian matrix. For
  any independent set $S$ of size $s$ and characteristic vector $z$, we have:
  \[ s \leq n \frac{\mu - \comp{d}_S}{\mu}.\]
  Furthermore, the following are equivalent:
  \begin{enumerate}[(a)]
    \setlength{\itemsep}{0pt}
    \setlength{\parskip}{0pt}
    \item Equality holds.
    \item $z$ is a linear combination of a $0$-eigenvector and a $\mu$-eigenvector.
    \item The bipartite subgraph induced by the partition $\{S,V(X)\setminus S\}$ is semi-regular.\qed
  \end{enumerate}
\end{theorem}

The $0$-eigenvectors of the Laplacian are exactly the constant vectors; compare this to \tref{Aboundeq}, where
the $k$-eigenvectors were the constant vectors, but only because the graph was regular. Note furthermore that
in both cases, when equality holds, $z$ is a linear combination of eigenvectors belonging to the greatest and
least eigenvalues. Based partly on the analogy between \tref{Aboundeq} and \tref{Lboundeq}, and the fact that
\tref{Aboundeq} is actually a special case of \tref{Lboundeq}, it seems that the Laplacian matrix formulation
is the natural generalization to non-regular graphs.

\vspace{2mm}

Assume that equality holds in \tref{Lboundeq} and $\gcd(s,n)=1$. The number of vertices in $S$ adjacent to a
vertex not in $S$, $\mu\frac{s}{n}$, must be an integer, and so $n \mid \mu$. As $0<\mu\leq n$, it follows
that $n=\mu$ and the bipartite subgraph induced by the partition $\{S,V(X)\setminus S\}$ is complete
bipartite. Of course the same conclusion follows if equality holds in \tref{Aboundeq}, since it is a special
case of \tref{Lboundeq}.

\section{Comparing Bounds}\label{sec:comp}

The bounds of \cref{Abound} and \cref{Lbound} are not directly comparable: we show this using
examples. We start with a simple family of graphs where \cref{Lbound} is better, then a family where
neither bound is uniformly better.

\vspace{2mm}

Consider the graphs $K_{a,b}$, where $a<b$. Clearly the only maximum independent set is the set of vertices of
degree $a$. So we have
\[ \comp{d}_S=a,\qquad k_S=\frac{2a^2}{a+b}. \]
Also, the least eigenvalue of the adjacency matrix is $-\sqrt{ab}$ and the greatest eigenvalue of the
Laplacian matrix is $a+b$. Applying \cref{Abound} we get
\begin{equation}\label{KabA}
  s \leq \frac{ (a+b)^2\sqrt{ab} }{ (a+b)\sqrt{ab} +2a^2}.
\end{equation}
However applying \cref{Lbound} we get exactly the size of the maximum independent set.
\begin{equation}
  s \leq b.\notag
\end{equation}
Thus we conclude that the bounds based on the adjacency and Laplacian matrices are not equal (for instance,
for $K_{4,23}$ \eqref{KabA} gives $s\leq24$). Furthermore, not only is \cref{Lbound} tight, but so is
\cref{Lbound2}. The latter bound is in terms of the graph only, whereas the former depends on $\comp{d}_S$ and
so retains an implicit dependence on the structure of $S$.

\vspace{2mm}

Now consider the graphs $X_m$, $m>1$, constructed as follows. Let $G_m$ be a copy of $\comp{K_m}$, and $H_m$ be a copy
of $C_{2m+1}$. Then $X_m$ consists of the disjoint union of $G_m$ and $H_m$, together with edges from every
vertex of $G_m$ to every vertex of $H_m$.

Clearly, the maximum independent sets are of size $s=m$ and there are two types: the vertices of
$G_m$ and an maximum independent set in $H_m$. The first type has $\comp{d}_S=2m+1$ and the second
has $\comp{d}_S=m+2$. Owing to the block structure of the adjacency and Laplacian matrices, the
eigenvectors can be determined. We find that $\tau = 1 - \sqrt{2m^2+m+1}$ and $\mu=3m+1$.

Using this information we can compute the two bounds. For convenience, let $\alpha_A$ and $\alpha_L$
be the values of the bounds in \cref{Abound} and \cref{Lbound}, respectively. If $S$ is the vertex
set of $G_m$; we determine that the Laplacian bound is tight while the adjacency bound is not.
\[ |S| = m = \alpha_L < \alpha_A. \]
\[ |S| = n\frac{\mu-\comp{d}_S}{\mu} < n\frac{-\tau}{k_S-\tau}. \]
Now let $S$ be a maximum independent set in $H_m$. We find that for small $m$, the adjacency bound is better,
and for large $m$ the Laplacian bound is better. Neither bound is ever tight, except that the adjacency bound
for $m=2$ is correct when rounded.
\[
  \def\arraystretch{2}
  \begin{array}{ll}
    |S|=m < \alpha_A < \alpha_L = \frac{2m-1}{3m+1} &\text{for }2 \leq m \leq 24 \\
    |S|=m < \alpha_L = \frac{2m-1}{3m+1} < \alpha_A &\text{for }25 \leq m
  \end{array}
\]

\subsection{Another Eigenvalue Bound}

We now compare the Delsarte-Hoffman bound of \cref{rbound} with another eigenvalue bound of Sarnak.

Let $\lambda$ be the maximum of the second largest eigenvalue and the absolute value of the least eigenvalue
of the adjacency matrix. Sarnak \cite{Sarnak} has shown the following bound for an independent set $S$ in a
$k$-regular graph.
\begin{lemma}\label{lem:sbound}
  \[ |S| \leq n \frac{\lambda}{k} \]
\end{lemma}
It turns out that this is strictly weaker than \cref{rbound}, as we now show. For reference, we give an
outline of the proof in \cite{Sarnak}.

\begin{proof}
  If $x \perp \1$ then $\norm{Ax}^2 \leq \norm{\lambda x}^2$. We choose
  \[ x_i = \begin{cases} n-s, &i \in S\\ -s, &i \notin S \end{cases} \]
  and compute the norms as follows:
  \begin{align}
  \norm{\lambda x}^2 &= \lambda^2 \pp{s(n-s)^2 + (n-s)s^2} = \lambda^2 ns(n-s), \notag\\
  \norm{Ax}^2 &= \sum_{i \in S} \pp{(Ax)_i}^2 + \sum_{i \notin S} \pp{(Ax)_i}^2 \phantom{\sum^X} \notag\\
    &\geq \sum_{i \in S} \pp{(Ax)_i}^2 \notag\\
    &= s^3k^2.\label{oops1}
  \end{align}
  Substituting into $\norm{Ax}^2 \leq \norm{\lambda x}^2$, we find that
  \[
  s \leq \sqrt{n(n-s)} \frac{\lambda}{k}
  \]
  which implies
  \begin{equation}\label{oops2}
    s \leq n \frac{\lambda}{k}.
  \end{equation}\qed
\end{proof}

Note that at \eqref{oops1} we are neglecting some positive terms, and at \eqref{oops2}, we are using $n-s \leq
n$. Assuming that $0<s<n$, either of these is sufficient to guarantee that the inequality in \lref{sbound} is
in fact strict. Based on these observations, we can improve this proof. If we write the adjacency matrix in
the form
\[ A = \mat{0&B\\B^T&C}, \]
then we can compute the missing contributions as follows:
\begin{align}
  \sum_{i \notin S} \pp{(Ax)_i}^2 &= \sum_i \pp{(n-s)(B^T\1)_i + (-s)(C\1)_i}^2 \notag\\
  &\geq \frac{1}{n-s} \pp{ \sum_i \pp{(n-s)(B^T\1)_i + (-s)(C\1)_i } }^2 \label{notoops}\\
  &= \frac{1}{n-s} \pp{ (n-s)(sk) + (-s)((n-2s)k) }^2 \notag\\
  &= \frac{s^4k^2}{n-s}. \notag
\end{align}
This gives that
\[ \norm{Ax} \geq \frac{s^3k^2n}{n-s}, \]
and then using $\norm{Ax}^2 \leq \norm{\lambda x}^2$ and rearranging gives exactly
\begin{equation}\label{sarnak2}
  s \leq \frac{n}{1+\frac{k}{\lambda}}.
\end{equation}
This bound is strictly better than \lref{sbound}, but, as $\lambda \geq -\tau$, it is no stronger than
\cref{rbound}. In the context of the present paper, this strengthening is not a surprise: the vector $x$ in
the proof of \lref{sbound} is just a multiple of $z-\frac{s}{n}\1$. Note also that if \eqref{sarnak2} holds
with equality then so does \eqref{notoops}: in the extremal case we haven't neglected any terms.

\section{A Generalization}

In developing the bounds of \sref{bounds}, we were motivated by a need to have tools we could apply to the
graphs $ER(q)$. It seems that our approach to the Delsarte-Hoffman bound can be pushed further. As one example
consider the following.

\begin{theorem}\label{thm:Bratio}
Let $X$ be any graph on $n$ vertices, and $S$ and independent set of size $s$. Let $B$ be a symmetric square
matrix indexed by the vertices of $X$ such that
\begin{enumerate}[(a)]
  \item $B \psd$
  \item $B_{ij} \leq 0$ whenever $i \neq j $ and $i \not\sim j$
  \item $B$ has constant row sum $r$
  \item $B$ has constant diagonal $t$
\end{enumerate}
Then
\[ s \leq n \frac{t}{r}. \]
\end{theorem}

This is relevant for the following reason. Let $\cA=\{\seq{A}{0}{d}\}$ be an association scheme and let $R
\sbs \{1,\ldots,d\}$ (see \cite[Chapter~12]{blue} for background and notation). An $R$-coclique is a set of
vertices such that no two of them are $i$-related for $i \in R$. In other words, it is an independent set in
the graph formed by the union of the classes of $R$.

We denote the trace of a matrix $N$ by $\tr(N)$, and the sum of all of its entries by $\summ(N)$. The
following result is shown in~\cite{asnotes}.

\begin{theorem}\label{thm:trsum}
Let $\cA$ be an association scheme with $d$ classes, let $R \sbs \{1,\ldots,d\}$ and let $S$ be an
$R$-coclique in $\cA$. Then
\[ |S| \leq \min \cc{ v \frac{\tr(N)}{\summ(N)} \mid N\psd, N \circ A_i \leq 0 \text{ for } i \notin R \cup
\{0\} } \]\qed
\end{theorem}

Let $X$ be the graph formed by the union of the classes of $R$. Then the matrix $N$ satisfies the criteria of
\tref{Bratio} with $tv=\tr(N)$ and $rv=\summ(N)$. So \tref{Bratio} is a generalization of \tref{trsum}.

On the other hand, we may write $N$ in terms of the matrix idempotents and Schur idempotents as
\[ N = \sum_i a_i A_i = \sum_j b_j E_j. \]
Thus $\tr(N)=va_0$ and $\summ(N)=vb_0$. Then \tref{trsum} says that $|S|$ is bounded above by the following
linear program.
\[ \min \cc{ v a_0 \mid b_j \geq 0, b_0=1, a_i \leq 0 \text{ for } i \notin R \cup \{0\} }. \]
This is equal to Delsarte's LP bound on an $\comp{R}$-coclique~\cite{Delsarte}. In other words, \tref{Bratio}
is a generalization of Delsarte's LP bound to general graphs.

\section{Erd\H{o}s-R\'{e}nyi}\label{sec:ER}

We now turn our attention to the motivation of our present work: applying the tools of \sref{bounds} to
derive new bounds on the independence number of the Erd\H{o}s-R\'{e}nyi graphs. This technique can in fact be
applied to more general polarity graphs.

\subsection{Quotients}\label{sec:quotients}

Let $Y$ be a $k$-regular graph, and let $\sigma$ be an automorphism of $Y$ or order $2$. The orbits of
$\sigma$ partition the vertices of $Y$. It can be seen by elementary means that this partition is
\emph{equitable}, that is, given any orbit $C_j$ and any vertex $u$, the number of edges from $u$ to a vertex
of $C_j$ depends only on the orbit $C_i$ containing $x$. Define $w_{ij}$ to be the number of edges between a
vertex of $C_i$ and the vertices of $C_j$. We define the \emph{quotient graph} $X = Y/\sigma$ to have vertex
set equal to the set of orbits of $\sigma$ with $w_{ij}$ arcs from $C_i$ to $C_j$.

Note that in general $X$ will be a weighted digraph, possibly with loops. It will be a graph (i.e.,
$w_{ij}=w_{ji}$) if and only if there is no edge between an orbit of size one and an orbit of size two. It
will have no multiple edges if and only if the edges joining two cells of size two never form a complete
bipartite graph. It will have no loops if and only if $\sigma$ never interchanges adjacent vertices.

The eigenvectors of $X$ correspond to the eigenvectors of $Y$ that are constant on each orbit of $\sigma$. In
particular, the eigenvalues of $X$ are exactly the eigenvalues of $Y$ corresponding to eigenvectors that are
constant on each cell (see \cite[Chapter~5]{blue} for details). Thus if we know the eigenvectors of $A(Y)$,
we know the least eigenvalue of $A=A(X)$, and we may apply \cref{Abound1}.

\subsection{Erd\H{o}s-R\'{e}nyi Graphs}

As a specific example, let $Y$ be the incidence graph of $PG(2,q)$. This is a bipartite graph, with points and
lines forming the two colour classes. Let $\sigma$ be the map that sends a point to the line with the same
coordinates, and vice versa. Then every orbit of $\sigma$ has size two and $Y / \sigma$ is an undirected graph
with no multiple edges. It does however have $q+1$ vertices with loops. The graph that results from removing
the loops is known as the \emph{Erd\H{o}s-R\'{e}nyi} graph of order $q$. For convenience, we will leave the
loops in. The graph $X$ has $q^2+q+1$ vertices, degree $q+1$, and eigenvalues $q+1,\pm\sqrt{q}$. Let $S$ be an
independent set of size $s$ containing $s_1$ loops in $X$. We use our results from \sref{bounds} to bound $s$.

\vspace{2mm}

We can apply \cref{rbound} directly; this will bound the size of an independent set containing no loops (i.e.,
an independent set that contains no absolute vertices). Our set $S$ could be at most $q+1$ greater than this,
yielding the following bound.
\begin{equation}\label{ratioboundERq}
  s \leq n \frac{-\tau}{k-\tau} + q+1 = \frac{\pp{q^2+q+1}\sqrt{q}}{q+\sqrt{q}+1} +q+1
\end{equation}
We can of course use \lref{sbound} instead of \cref{rbound}; this gives a manifestly weaker bound on $s$.

\vspace{2mm}

A better approach is to use \cref{Abound1}. In order to obtain a bound independent of $s_1$, we set
$s_1=q+1$, to get the following.
\begin{equation}\label{Abound1ERq}
  s \leq \frac{ \sqrt{q} + \sqrt{q + 4(q+1)\frac{q+\sqrt{q}+1}{q^2+q+1}} }{ 2 \frac{q+\sqrt{q}+1}{q^2+q+1} }
\end{equation}
As noted in \sref{Lapbounds}, this is equivalent to using \cref{Lbound} (or more precisely,
\cref{Lbound2}). Some tedious algebra shows that \eqref{Abound1ERq} is strictly better than
\eqref{ratioboundERq}.

\vspace{2mm}

Another approach would be to consider the graph $X_0$, obtained by deleting the absolute vertices from
$X$. Godsil and Royle have computed the characteristic polynomial of the graph $X_0$ to be
\[ (\lambda-q) \lambda (\lambda+1)^q (\lambda^2-q)^{(q^2-q-2)/2}. \]
So the least eigenvalue of $X_0$ is $-\sqrt{q}$. Let $S_0$ be an independent set of average degree
$\comp{d}_{S_0}$ in $X_0$. Trivially, we see that
\[ \comp{d}_{S_0} \geq q-1, \]
and so we can apply \cref{Abound} to $X_0$ and add $q+1$ to obtain the following bound.
\begin{equation}\label{noabsbound1ERq}
  s \leq \frac { q^2\sqrt{q} } { q-2+\frac{1}{q} + \sqrt{q} } + q+1
\end{equation}

On the other hand, we may assume that in fact $S_0$ in contained in an independent set $S$ of $X$. Thus $S$
consists of $S_0$ together with, say, $s_1$ absolute vertices. Each absolute vertex in $X$ is adjacent with at
most $q$ vertices of $S_0$, so there are at most $q(q+1-s_1)$ edges between $S_0$ and the set of absolute
vertices of $X$. It follows that
\[ \comp{d}_{S_0} \geq q+1 - \frac { q(q+1-s_1) } { s-s_1 }. \]
Substituting this into \cref{Abound} and adding $q+1$, we get a bound in terms of $s_1$. For $q \geq 5$, this is a
decreasing function of $s_1$, so we set $s_1=0$ to obtain the following.
\begin{equation}\label{noabsbound2ERq}
  s \leq \frac {q^2\sqrt{q}+2q(q+1)} {q+2+\frac{1}{q}+\sqrt{q}}
\end{equation}
For $5 \leq q \leq 23$, \eqref{noabsbound2ERq} is better than \eqref{noabsbound1ERq}, but for $q \geq 25$, the
reverse is true. Neither bound is as good as \eqref{Abound1ERq}.

\vspace{2mm}

Yet another approach is to delete the loop-edges from $X$ and apply \cref{Abound}. (This is the more usual
Erd\H{o}s-R\'{e}nyi graph.) Godsil and Royle have computed its characteristic polynomial to
be
\[ (\lambda^3-q\lambda^2-2q\lambda+q^2+q) (\lambda^2+\lambda+1-q)^q (\lambda^2-q)^{(q^2-q-2)/2}. \]
The least eigenvalue is a root of the cubic factor. We can approximate it using Newton's method. It is less
than $-\sqrt{q}$, it is the only eigenvalue that is less than $-\sqrt{q}$, and the cubic factor is concave
down for $\lambda \leq -\sqrt{q}$. So iterating Newton's method starting with $-\sqrt{q}$ will always give a
lower bound on the least eigenvalue, which means we will be overestimating our upper bound on the size of an
independent set. In fact, since we only care about the integer part of the final answer, it seems that two
iterations is sufficient.

Let $S$ be an independent set of size $s$ containing $s_1$ absolute vertices (here these no longer have
loops). It is straightforward to compute $\comp{d}_S$:
\[ \comp{d}_S = q + 1 - \frac{s_1}{s}. \]
Letting $w$ be an approximation to the least eigenvalue and applying \cref{Abound} we get
\begin{equation}\label{noloopboundERq}
s \leq \frac {(q^2+q+1)(-w) + q+1} {q+1-w +\frac{q+1}{q^2+q+1}}.
\end{equation}

We close this section with a brief table summarizing the numerical values of the bounds we have derived. We
also include exact values for the size of a maximum independent set; these are from
Williford~\cite[Section~4.3]{Williford}. In fact, it is partly the difference between his values and the bound
of \eqref{ratioboundERq} that motivated our work. The best bound we know of is \eqref{Abound1ERq}.

\begin{center}
\begin{tabular}{|cc|ccccc|}
\hline
$q$ & $\alpha(ER(q))$ & \eqref{Abound1ERq} & \eqref{noloopboundERq} & \eqref{ratioboundERq} &
        \eqref{noabsbound1ERq} & \eqref{noabsbound2ERq} \\
\hline
\hline
 3 &  5 &  5.56 &  5.63 &  7.92 &  9.09 &  6.21 \\
 5 & 10 & 10.56 & 10.82 & 14.42 & 16.28 & 12.28 \\
 7 & 15 & 16.73 & 17.27 & 22.16 & 24.65 & 20.50 \\
 9 & 22 & 23.93 & 24.87 & 31    & 34.03 & 29.98 \\
11 & 29 & 32.05 & 33.40 & 40.79 & 44.34 & 40.55 \\
13 & 38 & 41.03 & 42.88 & 51.48 & 55.49 & 52.08 \\
\hline
\end{tabular}
\end{center}

\subsection{Polarity Graphs}

Much of the work in the previous section can in fact be applied more generally.

Let $Y$ be the incidence graph of $PG(2,q)$, let $\sigma$ be an automorphism of $Y$ of order two that swaps
points and lines of $PG(2,q)$, and let $X = Y / \sigma$. The vertices in $X$ with loops are the absolute
vertices; denote their number by $a$. Let $N$ be the adjacency matrix of $X$. Then $N$ is an incidence matrix
of $PG(2,q)$ (the image under $\sigma$ of the point corresponding to the $i$-th row of $N$ is the line
corresponding to the $i$-th column of $N$), and the adjacency matrix of $Y$ is
\[ \mat{0&N\\N&0}. \]
In $PG(2,q)$ any point lies on $q+1$ lines and any two points lie on exactly one common line. It follows that
\[ N^2 = (q+1)I + (J-I). \]
Thus the eigenvalues of $N^2$ are $(q+1)^2$ and $q$ with multiplicities $1$ and $q^2+q$, and the eigenvalues
of $N$ are $q+1$, $\sqrt{q}$ and $-\sqrt{q}$ with multiplicities $1$, $m_1$ and $m_2$. The values of $m_1$ and
$m_2$ depend on $a$, but they are both non-zero. Thus we can apply \cref{Abound1} to bound the size of an
independent set $S$ in $X$.
\[ |S| \leq  \frac{ \sqrt{q} + \sqrt{q + 4a\frac{q+\sqrt{q}+1}{q^2+q+1}} }{ 2 \frac{q+\sqrt{q}+1}{q^2+q+1} } \]

\vspace{2mm}

A similar approach yields bounds on the size of an independent set in polarities of the generalized
quadrangles $W(q)$. This consists of a set of points, a set of lines, and an incidence relation between them
such that
\begin{enumerate}[(a)]
  \item Each point is incident with $q+1$ lines and two distinct points are incident with at most one common line.
  \item Each line is incident with $q+1$ points and two distinct lines are incident with at most one common point.
    \item Given a point $p$ and a line $l$ not incident with $p$, there exists a unique point $q$ and a unique
    line $m$ such that $m$ is incident with $p$ and $q$ and $q$ is incident with $m$ and $l$.
\end{enumerate}
The reader is directed to \cite{PT} for more details. The number of points in $W(q)$ is $q^3+q^2+q+1$, which
is also the number of lines. Thus a polarity graph of $W(q)$ has $q^3+q^2+q+1$ vertices and degree $q+1$
(leaving the loops in). Computing the eigenvalues of a polarity graph of $W(q)$ is a little more work than for
a polarity graph of $PG(2,q)$, but the argument is similar (see \cite[Section~1.8.2]{PT} for details). The
least eigenvalue is $-\sqrt{2q}$. Furthermore, the number of absolute points is always $q^2+1$. Thus
\cref{Abound1} gives the following bound on the size of an independent set $S$ in any polarity graph of
$W(q)$.
\[
|S| \leq \frac { \sqrt{2q} + \sqrt{ 2q + 4(q^2+1)\frac{q+\sqrt{2q}+1}{q^3+q^2+q+1} } }
{ 2 \frac{q+\sqrt{2q}+1}{q^3+q^2+q+1} }.
\]

\bibliographystyle{plain}
\bibliography{ratio}

\end{document}